\magnification 1200
\def\R{{\rm I\kern-0.2em R\kern0.2em \kern-0.2em}}
\def\N{{\rm I\kern-0.2em N\kern0.2em \kern-0.2em}}
\def\P{{\rm I\kern-0.2em P\kern0.2em \kern-0.2em}}
\def\B{{\rm I\kern-0.2em B\kern0.2em \kern-0.2em}}
\def\Z{{\rm I\kern-0.2em Z\kern0.2em \kern-0.2em}}
\def\C{{\bf \rm C}\kern-.4em {\vrule height1.4ex width.08em depth-.04ex}\;}

\def\z{{\zeta}}

\def\cH{{\cal H}}

\font\ninerm=cmr8
\noindent  {\ninerm To appear in Proc.\ Roy.\ Soc.\ Edinb.\ Sect.\ A} 
\vskip 25mm

\centerline {\bf SINGLE VALUED CONJUGATES AND}

\centerline {\bf HOLOMORPHIC EXTENDIBILITY}
\vskip 4mm
\centerline{Josip Globevnik}
\vskip 4mm
{\noindent \ninerm ABSTRACT\ \  Let $D$ be a bounded domain in the complex plane 
whose boundary consists of $m\geq 2$ pairwise disjoint simple closed curves and let   
$A(bD)$ be the algebra of all continuous 
functions on $bD$ which extend holomorphically through $D$. We show that a 
continuous function $\Phi $ on 
$bD$ belongs to $A(bD)$ if for each $g\in A(bD)$ 
the harmonic extension of $\Re (g\Phi )$ to $D$ has a single valued conjugate.
} 
\vskip 6mm
Let $D$ be a bounded domain in the complex plane whose boundary consists of $m$ pairwise 
disjoint simple closed curves.  A continuous real function $\varphi $ on $bD$ 
has a unique continuous extension $\cH (\varphi )$ to $\overline  D$ which is harmonic on $D$. 
The conjugate $\cH (\varphi )^\ast $ of $\cH (\varphi )$ is a harmonic function on $D$ which,
in general, is multiple valued. It is determined up to an additive constant. 
Let $A(bD)$ be the algebra of all continuous functions 
on $bD$ which extend holomorphically through $D$.

If $\Phi \in A(bD)$ then the harmonic extension of $\Re\Phi$, the real part of $\Phi $, 
has the single valued conjugate
$\Im\tilde\Phi $ where $\tilde\Phi $ is the holomorphic extension of $\Phi $ through $D$. 
Thus, for any $g\in A(bD)$, the harmonic extension of $\Re (g\Phi )$ has
a single valued conjugate.  In the present note we show that in the case when 
$D$ is multiply connected, 
that is, when $m\geq 2$ this property characterizes the functions from $A(bD)$. 
\vskip 2mm
\noindent\bf PROPOSITION \it \  Let $D$ be a bounded domain in the plane whose
boundary consists 
of $m\geq 2$ pairwise disjoint simple closed curves. A continuous function $\Phi $ on $bD$ 
extends holomorphically through $D$ if and only if for each $g\in A(bD)$  the  
conjugate of the harmonic extension of  $\Re (g\Phi )$ is single valued. \rm
\vskip 2mm
\noindent In other words, $\Phi $ extends holomorphically through $D$ if and only 
if for each $g\in A(bD)$ the harmonic extension of $\Re (g\Phi)$ to $D$ is the real 
part of a single valued holomorphic function on $D$.

There is no such proposition in the case when $D$ is simply connected as in this case 
every harmonic function on $D$ has a single valued conjugate.
  
Every bounded domain $D\subset \C$ whose boundary consists of finitely 
many pairwise disjoint simple closed 
curves is biholomorphically equivalent to a domain $D^\prime$ bounded by finitely many 
pairwise disjoint circles [4]. Moreover, every biholomorphic map 
$\Phi\colon D\rightarrow D^\prime $ extends to 
a homeomorphism $\tilde\Phi\colon \overline D\rightarrow\overline{D^\prime}$ 
(see the proof in [2,\ pp.\ 46-49] which works 
also for multiply connected domains, bounded by finitely many pairwise disjoint 
simple closed curves). 
Thus, with no loss 
of generality assume that $D$ is bounded by finitely many pairwise disjoint circles  
$\gamma _1, \gamma _2,\cdots, 
\gamma _m, \ m\geq 2$ where $\gamma _m$ bounds the unbounded component 
of $\C\setminus\overline D$.  For each $j,\ 1\leq j\leq m,$ let $h_j$ 
be the continuous function on $\overline D$, harmonic on $D$, such that
$h_j\vert \gamma_j\equiv 1$,\ $h_j\vert\gamma _k\equiv 0\ (1\leq k\leq m,\ k\not = j)$. These
functions  extend harmonically across $bD$. For each $j,\ 
1\leq j\leq m$ let $\gamma ^\prime_j 
\subset D $ be a circle that has the same center as $\gamma _j$ and that is a slight 
perturbation of $\gamma_j$. If $u$ is a real harmonic function on a neighbourhood of 
$\overline D$ then the period of 
$u^\ast $ along  $\gamma ^\prime_j$ is the same
as the period of $u^\ast $ along $\gamma _j$, that is, equal to 
$\int _{bD}\varphi {{\partial h_j}\over{\partial n}}
ds $ [3, p.80]. Since every real continuous function on $\overline D$ which is harmonic on $D$ is
the uniform limit of a sequence $u_n$ such that 
each $u_n$ is harmonic on a neighbourhood of $\overline D$ [3, p. 79] the preceding dicussion 
implies that if $\varphi $ is a real 
valued continuous 
function on $bD$ then the period of $\cH (\varphi)^\ast$ along $\gamma _j^\prime$ equals 
$\int _{bD}\varphi {{\partial h_j}\over{\partial n}}
ds $ so $\cH (\varphi )^\ast $ is 
single valued if and only if
$$
\int _{bD}\varphi {{\partial h_j}\over{\partial n}}ds = 0\ \ (1\leq j\leq m-1). 
\eqno (1)
$$
For each $j, \ 1\leq j\leq m$, \ $h_j+ih_j^\ast $ is a multiple 
valued holomorphic function in a neighbourhood of $\overline D$ whose complex derivative $W_j=
(h_j+ih_j^\ast)^\prime $ is a single valued holomorphic function in a neighbourhood of 
$\overline D$ without 
zeros on $bD$ [3, p.81]. The Cauchy-Riemann equations imply that on $bD$
$$
{{\partial h_j}\over{\partial n}}ds = iW_j(\z )d\z \ \ (1\leq j\leq m)
$$
[3, p.81]. Hence, by (1), $\cH (\varphi )^\ast $ is single valued 
if and only if 
$$
\int _{bD}\varphi (\z )iW_j (\z )d\z = 0\ \ (1\leq j\leq m-1).
\eqno (2)
$$

The only if part of the proposition follows from the discussion preceding the proposition.
To prove the if part, assume that $\Phi $ is a continuous function on $bD$ such that  
for each $g\in A(bD)$  the  
conjugate of the harmonic extension of  $\Re (g\Phi )$ is single valued. 

By (2) it follows that for each $g\in A(bD)$ we have
$$
\int_{bD} \Re [g(\z )\Phi (\z )]iW_j(\z )d\z = 0\ \ (1\leq j\leq m-1)
$$
which, since $iW_j(\z )d\z $ is real on $bD,\ 1\leq j\leq m-1$, implies that
$$
\Re\Bigl[\int_{bD}g(\z )\Phi (\z )iW_j (
\z )d\z \Bigr] = 0\ \ (1\leq j\leq m-1).
$$
It follows that for each $z\in \C\setminus \overline D$ we have
$$
\int _{bD} {{1}\over{\z - z}}\Phi (\z )W_j(\z )d\z = 0\ \ (1\leq j\leq m-1)
$$
which implies that for each $j,\ 1\leq j\leq m-1,$\ the function $\z\mapsto \Phi (\z)W_j (\z )$ 
extends holomorphically through $D$. So, for each $j,\ 1\leq j\leq m-1$, there is a 
function $H_j\in A(bD)$ such that 
$$
\Phi (\z ) = {{H_j(\z )}\over{W_j(\z )}}\ \ (\z \in bD,\ 1\leq j\leq m-1)
\eqno (3)
$$
which means that for each $j,\ 1\leq j\leq m-1$, the function $\Phi $ extends 
through $D$ as a meromorphic function with possible poles at the zeros of the function 
$W_j$. If we show that $W_j,\ 1\leq j\leq m-1$, have no common zero on D then (3) will
imply that $\Phi $ extends holomorphically through $D$. 

Fix $a\in D$ so close to $bD$ that all the zeros $a_1, a_2,\cdots a_{m-1}$ of the Szeg\" o 
kernel $z\mapsto S(z,a)$ are simple [1, p.106]. If $w\in D$ then the Garabedian kernel 
[1, p.24] is a meromorphic function on D with a single simple pole at $w$ [1, p.25]
which has no zero on $D\setminus\{ w \}$ [1, p.49]. 
We now use the fact that the complex linear span of $\{W_1, W_2,\cdots W_{m-1}\} $ 
coincides with the complex linear span of\ 
$
\{ L(z,a_j)S(z,a):\ 1\leq j\leq m-1\}
$ 
[1, p.80]. The preceding discussion implies that for each $j,\ 1\leq j\leq m-1$, 
the function  $L(z, a_j)S(z,a)$ has zeros precisely at $a_k,\ 1\leq k\leq m-1,\ k\not= j$,
which 
implies that the functions $L(z, a_j)S(z,a), 1\leq j\leq m-1$ have no common zero on $D$. 
Suppose that
$W_j(b)=0\ (1\leq j\leq m-1)$ for some $b\in D$. Then every function in the complex linear 
span of $\{W_1, W_2,\cdots W_{m-1}\} $ vanishes at $b$. In particular, all functions
 $L(z, a_j)S(z,a), 1\leq j\leq m-1$, vanish at $b$, a contradiction. This shows that 
 the functions 
$W_j,\ 1\leq j\leq m-1$, have no common zero on D. The proof of the proposition is complete.  
\vskip 2mm
We conclude with an example. Let $0<R<1 $ and let $D= \{ \z\colon \ R<|\z |<1\}$. Let 
$\Phi = \varphi |bD$  where 
$\varphi (z)=\overline z$. Let $g(z)=z^n,\ n\in Z,\ n\not= 1$. Then $\Re (g\Phi) (e^{i \theta }) =
\cos ((n-1)\theta) $,\ \ $
\Re (g\Phi) (Re^{i \theta }) = R^{n+1}\cos ((n-1)\theta )\ \ \  (0\leq \theta<2\pi )$. The function
$$
u(z) = 
\Re\Bigl[{{R^{n+1}-R^{n-1}}\over{R^{n-1}-R^{1-n}}}(z^{n-1}-z^{1-n})+ z^{n-1}\Bigr]
$$
is harmonic on $\C\setminus\{ 0\}$ and satisfies $u(e^{i\theta }) = \cos ((n-1)\theta) $, \ \ 
$u(Re^{i \theta }) =  R^{n+1}\cos ((n-1)\theta)\ \ \   (0\leq \theta<2\pi )$\  
so $u$ provides a harmonic extension of the function 
$\z\mapsto \Re [\z ^n\Phi (\z )] \ (\z\in bD)$ to $\C\setminus\{ 0\}$. Clearly 
$$
u^\ast (z) = \Im\Bigl[{{R^{n+1}-R^{n-1}}\over{R^{n-1}-R^{1-n}}}(z^{n-1}-z^{1-n})+
z^{n-1}\Bigr]
$$ 
is single valued. Thus, for every $n\in Z,\ n\not= 1$, the conjugate of the harmonic extension of the function 
$
\z\mapsto \Re [\z ^n\Phi (\z )] \ (\z\in bD)$ is single valued on $D$, yet 
$\Phi $ does not extend holomorphically through $D$. 
\vskip 4mm
This work was supported 
in part by the Ministry of Higher Education, Science and Technology of Slovenia  
through the research program Analysis and Geometry, Contract No.\ P1-0291. 
\vfill
\eject
\centerline{\bf REFERENCES} \rm
\vskip 6mm
\noindent 
1. S.\ Bell:\ \it The Cauchy Transform, Potential Theory and Conformal Mapping. \rm  
CRC Press, 
Boca Raton - Ann Arbor - London - Tokyo 1992 
\vskip 2mm
\noindent 
2. E.\ F.\ Collingwood and A.\ J.\ Lohwater: \it The Theory of Cluster Sets. \rm
 Cambridge University Press, Cambridge 1966
\vskip 2mm
\noindent
3. S.\ D.\ Fisher:\ \it Function Theory on Planar Domains.\rm  John Wiley \& Sons, New York 1983
\vskip 2mm
\noindent 
4. G.\ M.\ Goluzin: \it Geometrische Funktionentheorie. \rm
\rm VEB Deutscher Verlag 
der Wissenschaften, Berlin 1957
\vskip 20mm
\noindent Institute of Mathematics, Physics and Mechanics

\noindent University of Ljubljana

\noindent Ljubljana, Slovenia

\noindent josip.globevnik@fmf.uni-lj.si

\bye